\documentclass[conference]{IEEEtran}
\usepackage{algorithmic}
\usepackage{algorithm}
\usepackage{mathrsfs,amsthm}
\usepackage{cite}
\IEEEoverridecommandlockouts
\makeatletter
\def\ps@headings{%
\def\@oddhead{\mbox{}\scriptsize\rightmark \hfil \thepage}%
\def\@evenhead{\scriptsize\thepage \hfil \leftmark\mbox{}}%
\def\@oddfoot{}%
\def\@evenfoot{}}
\makeatother 
\usepackage{graphicx,amssymb,stfloats,subfigure,multicol}
\usepackage{epsfig,epsf,psfrag,amssymb,amsfonts,latexsym,graphicx,mathrsfs,subfigure}
\usepackage{bbm}
\usepackage{chngpage}
\usepackage[dvips]{color}
\usepackage{textcomp}
\usepackage{hyperref}
\usepackage[hyphenbreaks]{breakurl}
\usepackage[normalem]{ulem}
\newtheorem{Def}{Definition}
\newtheorem{corollary}{Corollary}
\newtheorem{thm}{Theorem}
\newtheorem{lem}{Lemma}
\newtheorem{prop}{Proposition}

 \def\old#1{}    

\def\beq{\begin{equation}}
\def\eeq{\end{equation}}
\def\bea{\begin{eqnarray}}
\def\eea{\end{eqnarray}}
\def\ba{\begin{array}}
\def\ea{\end{array}}

\def\bitem{\begin{itemize}}
\def\eitem{\end{itemize}}
\def\ben{\begin{enumerate}}
\def\een{\end{enumerate}}

\def\eg{{\it e.g., \/}}

\def\ie{{\it i.e.,\ \/}}

\definecolor{bgrd}{rgb}{1,1,1}
\definecolor{gray}{rgb}{0.5,0.5,0.5}
\definecolor{dkr}{rgb}{0.7,0.1,0.2}
\definecolor{dkb}{rgb}{0.1,0.1,0.8}

\def\edoc{\end{document}}

\newcommand{\mbbE}{\mathbb{E}}

\newcommand{\Tmsc}{\mathscr{T}}

\begin{document}

\title{On the Dynamics of Distributed Energy  Adoption:
Equilibrium, Stability, and Limiting Capacity }
\author{\large Tao Sun, Lang Tong, and Donghan Feng
\thanks{\scriptsize
 Tao Sun is with the Department of Civil and Environmental Engineering, Stanford University, Stanford, CA 94305, USA. }
\thanks{\scriptsize Lang Tong is with the School of Electrical and Computer
Engineering, Cornell University, Ithaca, NY 14853, USA Email: lt35@cornell.edu.}
\thanks{\scriptsize
 Donghan Feng is with the Department of Electrical Engineering, Shanghai Jiao Tong University, Shanghai 200240, P.R. China. }
\thanks{\scriptsize This work is supported in part by the National Science Foundation under Awards
	1809830 and 1816397.}
\thanks{A preliminary version this paper was published as a conference paper in \cite{Sun&Tong:17Allerton}}}

\maketitle

\begin{abstract}
The death spiral hypothesis in electric utility represents a positive feedback phenomenon in which a regulated utility is driven to financial instability by rising prices and declining demand.   We establish conditions for the existence of death spiral and conditions of stable adoption of  distributed energy resources.   We show in particular that linear tariffs always induce death spiral when the fixed operating cost of the utility rises beyond a certain threshold.
 For two-part tariffs with connection and volumetric charges, the Ramsey pricing that optimizes myopically social welfare subject to the revenue adequacy constraint  induces a stable equilibrium. The Ramsey pricing, however, inhibits renewable adoption with a high connection charge.  In contrast, a two-part tariff with a small connection charge results in a stable adoption process with a higher level  of renewable adoption and greater long-term total consumer surplus.  Market data are used to illustrate various solar adoption scenarios.

\end{abstract}


\section{Introduction}
Death spiral for electric utilities stands for a positive feedback scenario  in which, when the utility raises its price to cover its cost, consumers reduce consumption.  This forces  the utility to increase further its price, which lowers the consumption even further.

The possibilities of death spiral for electric utilities have been raised several times since 1960's \cite{Costello&Hemphill:14EJ}, and this topic has attracted considerable attention recently, thanks to the rapid deployment of the behind-the-meter solar photovoltaic (PV) and other distributed energy resources such as storage.   A main difference this time is the role of disruptive technology such as solar PV and residential storage.  Both technologies have direct impacts on the revenue of the utility.

There is some evidence supporting the underlying assumptions of the death spiral hypothesis.  Recent reports issued by the California Public Utility Commission (CPUC) \cite{CPUC:17_1,CPUC:17_2} state that ``From 2012 to 2016, system average rates (SAR) across the three IOUs has increased at an annual average of approximately 3.44\%, which is well above the average annual inflation rate of 1.3\% over the same time period.''  Meanwhile, ``all three utilities have experienced declines in kWh sales, which also lead to increased rates when revenue requirement remains flat or rises.''  Data in \cite{GosolarCA:17} further show  that  ``the
flattening or declining trend in kWh sales is driven by a changing
economy, growth in the customer (so called behind-the-meter) solar
industry, increasing availability of demand side management
programs such as energy efficiency, and the incremental proliferation of
retail choice.''

The above snapshot statistics  are consistent with the more general trend discussed in one of the earliest work on death spiral hypothesis in solar PV adoption by Cai, Adlakha, Low, Martini, and Chandy \cite{CaiEtal:13EP}.  Using data from an investor-owned regulated utility, the results in  \cite{CaiEtal:13EP} show, empirically, effects of positive feedback loop on  PV adoption, the loss of revenue,
and rate changes.  The empirical analysis also shows  that high connection charges slow the rate of solar adoption.
A more recent empirical study \cite{DarghouthEtal:16AE} using nation-wide data  by Darghouth, Wiser, Barbose, and Mills, besides confirming the general feedback phenomenon and the negative impact of connection charges on PV adoption,  shows more nuanced effects of dynamic pricing  on PV adoption.

While empirical studies suggest the potential of death spiral,  they lack the predictive power on the dynamics of the feedback loop of PV adoption and its policy implications.  With decreasing costs of solar PV, there is a need for a fundamental understanding of the PV adoption dynamics and impacts of key parameters in the adoption process.  Such parameters include the cost of solar, tax incentives, and the operating cost of the utility.

\subsection{Summary of Results}  This paper complements existing empirical studies such as \cite{CaiEtal:13EP,DarghouthEtal:16AE} with an analytical study on the dynamics of PV adoption.  In particular, we aim to shed lights on the following questions:
\bitem
\item Can death spiral happen under the current tariff?
\item What are the conditions and pricing mechanisms for a stable adoption\footnote{In economics, the process of adoption of technology is also often referred to technology {\em diffusion} \cite{Rogers:03book}.  In this paper, we use the word adoption in place of diffusion to avoid a conflict of terminology with the control theoretic notion of diffusion.} of renewable technology?
\item What is the maximum installation capacity (referred to as the limiting capacity) achievable by a stable adoption process?
\item Does a higher level of renewable penetration imply greater social welfare?
\eitem
The main contribution of this work is an analytical framework that allows us to study the PV adoption process
as a nonlinear dynamical system.  This model captures interactions between a regulated utility and its price-elastic and  net-metered\footnote{Net-metering stands for the mechanism where a consumer is charged by the net-consumption of electricity, which is equivalent to allowing a consumer sell locally generated electricity back to the utility.} consumers
who maximize consumer surplus and  make PV adoption decisions based on the payback time of the solar investments.  Such decisions are influenced by the tariff set by the regulator and the cost of solar PV.

By analyzing the nonlinear system with the tariff and the installed solar capacity as its states, we establish four main theorems on the dynamics of the solar PV adoption process: (i) the existence condition of death spiral in Theorem \ref{th:deathcond} and its application to the Ramsey linear tariff; (ii) the condition for stable adoption in Theorem \ref{th:stability}; (iii) the stability of the adoption process induced by the Ramsey two-part tariff in Theorem \ref{th:ramseytwopart}; and (iv)  the achievability of limiting adoption capacity in Theorem \ref{deflim}. These conditions are applied to benchmark tariff policies.  

A main conclusion of this work is that linear tariffs are prone to induce death spiral; so are the two-part tariffs with arbitrarily set connection charges.  On the other hand, the Ramsey pricing with the optimized volumetric and connection charges  guarantees a stable adoption. The high connection charge of the Ramsey price, however,  has a negative impact on PV adoption. We show, in fact, that  Ramsey pricing stalls PV adoption. In contrast, a mechanism that adds a small connection charge induces a stable adoption process that achieves a higher level of PV adoption.  We demonstrate in addition that, while  maximizing the immediate overall consumer surplus, Ramsey pricing  may generate smaller consumer surplus in the long run.

We also report an empirical study using wholesale and retail prices, demand data, estimated revenue, consumption profiles, and roof-top solar in New York city.  In this setting, we study the potential of death spiral and effects of tariff on solar adoption and consumer surplus under short-run and long-run analysis models.  The main conclusion of this study is twofold.  First, the default tariff  (a two-part tariff)  by Consolidated Edison (ConEd) of New York does not lead to death spiral, and it offers a comfortable stability margin. Second, there is a potential to increase the level of stable adoption and consumer surplus in the long run if the current connection charge is lowered judiciously. 

The literature is limited on the {\em dynamics} of PV adoption  beyond the empirical studies in \cite{CaiEtal:13EP,DarghouthEtal:16AE} and economic analysis \cite{Costello&Hemphill:14EJ}.  To our best knowledge,  this work appears to be the first to pursue an analytical characterization of the PV adoption dynamics in the framework of nonlinear dynamical feedback systems.

In many ways, whereas our results corroborate conclusions in  \cite{CaiEtal:13EP,Costello&Hemphill:14EJ,DarghouthEtal:16AE}, we provide deeper analytical insights into the role of tariff on the adoption process including how connection charges affect the level of adoption and  ways to mitigate the threat of death spiral. For instance, there is a consensus that, although the possibility of death spiral is real in the era of greater DER, the likelihood of a death spiral occurring is small, especially if the regulator and the utility set the tariff policy proactively, including the proper use of connection charges \cite{Costello&Hemphill:14EJ}. Our analysis is consistent with these conclusions. We provide, however, qualitative and  quantitative answers on how such proactive measures can be applied dynamically.  

\subsection{Related Work}

The design of retail tariff in the distribution system is an instance of the classical pricing problem for a regulated monopoly \cite{Brown&Sibley:86book}. In approving a proposed tariff, the regulator takes into account the impact of the tariff on overall social welfare, fairness, and societal concerns.  In such a setting,  the classical Ramsey pricing aims to maximize the social welfare subject to the break-even constraint for the utility \cite{Ramsey:27EJ,Brown&Sibley:86book}.  In this context, we consider the class of linear tariffs and the class of (nonlinear) two-part tariffs defined by a connection charge and  a vector volumetric charge.   Originally studied by Oi in his seminal work \cite{Oi:71QJE}, the two-part tariff is now widely adopted by utilities for residential customers in the United States  where  nearly 87\% of the residential customers face some form of connection charges \cite{Borenstein:16EJ}.

Tariff models for electricity markets with stochastic demand  have been extensively studied. See \eg \cite{Joskow&Tirole:06RJE} and references therein.  With the increasing presence of distributed energy resources (DER), there is heightened attention on different types of tariffs \cite{Costello:15EJ}.  In such settings, the Ramsey pricing problem for the retail utility in distribution systems with stochastic distributed energy resources is considered in \cite{Jia&Tong:16TSG,Jia&Tong:16JSAC,Munoz-Alvarez&Tong:17TPS1,Munoz-Alvarez&Tong:17TPS2}.  Our dynamic model builds upon the analysis of retail tariff design in  \cite{Munoz-Alvarez&Tong:17TPS1,Munoz-Alvarez&Tong:17TPS2}.

A key component of our analysis is to incorporate a solar PV adoption  model. Prior studies have modeled the PV adoption in two ways. A number of them are based directly on the discrete choice model\cite{MaribuEtal:07EP,CaiEtal:13EP}. Other studies apply existing adoption models for innovation that capture higher level characteristics of the so-called $S$-curve \cite{Bass:69MS} such as market potential and adoption fraction\cite{NREL:09,Islam:14EP,DarghouthEtal:16AE}. The separable formulation of market potential in the second way fits directly to the stability analysis of the investigated dynamical system. To this end, we adopt a widely used $S$-curve model for the aggregated consumer behavior \cite{Rogers:03book,Link:99book,Mahajan&Muller&Bass:91bookchap}, under an implicit assumption of successful PV adoption.

\section{Consumer, Retailer, and Adoption Models} \label{sec:II}

\subsection{Retail Tariff Structure}
We consider  retail tariffs that are uniformly applied to all consumers.  We assume that the retailer sets tariff $T$ that is subject to approval by the regulator periodically, say, on a daily, monthly or  yearly basis.  In the rate-setting period $k$, the tariff $T_k$ is fixed until the next period.

Mathematically,   {\em tariff} is a {\em pricing policy} that maps consumptions to payment.  To this end, we consider two classes of tariffs that are widely used in practice:

Two widely applied tariff classes are considered:

\ben
\item {\em Linear tariff :}   $\Tmsc_{\rm L}=\{T: T(d)=\pi^\top d\}$   where $d\in \mathbb{R}^N$ is a vector of consumptions. Here $N$ is the number of consumption periods in a billing cycle and d the vector of consumptions referred to as the (load) consumption profile.  Vector $\pi\in \mathbb{R}^N$ is the  the vector of marginal prices of electricity.   We say price $\pi$ is a time varying (or dynamic) when the consumptions are priced differently over time, and the price is {\em flat} when all entries of $\pi$ are the same.       The class of {\em flat tariffs} is thus defined as $\Tmsc_{\rm F}=\{T: T(d)=(\pi_{\rm F} {\bf 1}^{\top} d\}$. 
\item {\em Two-part (affine) tariff:} $\Tmsc_{\rm A}=\{T: T(d)=A+\pi^\top d\}$ where $A$ is the connection charge independent of the consumption.

\een
Evidently, $\Tmsc_{\rm F}$ is a subclass of $\Tmsc_{\rm L}$, and $\Tmsc_{\rm L}$ a subclass of  $\Tmsc_{\rm A}$.

\subsection{Consumer Decision Model}

We assume price-elastic demands, and consumer $i$'s demand depends on the local random state $\omega_i\in \mathbb{R}^N$ that is assumed to be an exogenous random process.

We assume that,  for a given tariff $T$ set by the retailer,   each  consumer decides its consumption by maximizing  the consumer surplus: it follows from \cite{Jia&Tong:16TSG,Munoz-Alvarez&Tong:17TPS1} that knowing the  tariff $T$, consumer $i$ who maximizes his surplus solves the multistage stochastic problem:
\begin{equation}
\label{eq:cs}
\overline{\rm cs}_i(T)=\mathop {\max }\limits_{q} \mathbb{E}\Bigg[u_i(q,\omega_i)-T(q-r_i(\omega_i))\Bigg],
\end{equation}
where  $u_i(q,\omega_i)$ is the utility of consuming $q$, and $r_i(\omega_i)$ the realized behind-the-meter renewable generation for consumer $i$.  In (1), $\overline{\rm cs}_i(T)$ is the (optimized) consumer surplus under tariff $T$, and the optimized consumption, denoted as $D_i(T,\omega_i)$, is consumer $i$'s load profile.  Note that the consumer's decision does not depend directly on the wholesale price electricity due to the implicit assumption that the consumer does not have access to real-time price in the wholesale market; once the price of electricity is set, the consumption decision depends only on the tariff and the local state variable.

With total $M$ consumers in the service area of the utility, the expected consumer surplus  under a two-part tariff  is
\bea
\label{eq:csaggre}
\overline {\rm cs} (T,R)
=\mathbb{E}[U(T,\omega) - \pi^\top(D(T,\omega)- R{r}_0(\omega))] - MA,
\eea 
where  $\omega=(\omega_1,\cdots, \omega_M)$ is the random state of all customers,
{$U(T,\omega) =\sum_i u_i(D_i(T,\omega_i),\omega_i)$ the aggregated utility and $D(T,\omega)=\sum_iD_i(T,\omega_i)$ the aggregated demand, respectively.     The first term on the right hand side of (\ref{eq:csaggre}) is the aggregated consumer utility, the second the total volumetric charge, and the last the total connection charge. The total renewable generated behind the meter is given by $R{r}_0(\omega) $ where  ${r}_0(\omega)$ is the renewable generation per unit-capacity installed  and $R$ the total installed capacity.

	\subsection{Retailer Decision Model} We model the retail utility as a regulated monopoly, which is the case in most parts of the United States.   Here we assume that
	the retailer imports electricity from the wholesale market to satisfy the aggregated demand of its customers.   The retailer is assumed to be a price taker\footnote{A large retail utility, strictly speaking, can influence the wholesale price of electricity.}.   This model is a reasonable approximation of the deregulated two-settlement electricity market.

	The retailer sets the tariff and seeks its approval by the regulator in each tariff-setting period.   As a regulated monopoly, the retailer is allowed to break even. The revenue adequacy condition is met by setting the {\em retail surplus} to zero. The retail surplus is defined as  		
	\beq
	\label{eq:rsexpression}
	\overline {\rm rs} (T,\theta,R) =\mbbE[ (\pi-\lambda)^\top(D(\pi ,\omega )-Rr_0(\omega))] +MA-\theta
\end{equation}
where $\theta$ is the operating cost of the utility, $\lambda\in R^N$ the wholesale price of electricity, and $(D(\pi ,\omega ))-Rr_0(\omega))$ the net consumption. The expectation is taken over $\lambda$ and $\omega$.
The first term on the right-hand side is the revenue from energy consumption.  The second term  $(MA)$ is the revenue from the connection charge.  The break-even condition can be satisfied by jointly allocating these two types of revenue to the operating cost of the utility.

Under the break-even constraint, the retailer has additional dimensions of  freedom to set the tariff to achieve a variety of objectives, including maximizing consumer surplus and enhancing the overall social welfare.  To this end, we model the retailer's pricing decision by a {\em tariff policy} $\mu$ that maps its expected future  operating cost $\theta$  and the current level of renewable adoption $R$ to a tariff $T$ in some tariff class in the next period.  In particular, at the end of the $k$th period, the tariff in the next period $T_{k+1}$ in either $\Tmsc_{\rm L}$ or $\Tmsc_{\rm A}$ is given by
\[
\mu:~~T_{k+1} = \mu(R_k,\theta_k)
\]
where  $R_k$ is the installed capacity at the end of period $k$ and $\theta_{k}$ is the utility's  expected fixed cost\footnote{Rigorously speaking, the tariff $T_{k+1}$ should take into account not only installed PV in $R_k$ but also new installations in period $k+1$, which depends on the cost of solar $\xi$.  For simplicity, we ignore the dependency of $\xi$.}.

An important type of tariff policy  is {\em Ramsey pricing} in which the retailer maximizes the total surplus\footnote{The sum of consumer and retail surpluses is sometimes referred to as social welfare.} subject to the revenue adequacy constraint.  Equivalently, the retailer solves the following constrained optimization to determine $T_{k+1}$  given
the current  level of renewable installation $R_k$ and the (expected) fixed cost $\theta_{k}$ in the next period:
\begin{equation}
\mu^*:~~ \mathop {\max }\limits_{T\in \mathscr{T}} \overline {\rm cs} (T,R_k)\quad
{\rm{s.t.}}\;\;\overline {\rm rs} (T,\theta_k,R_k) = 0.
\label{eq:ramplan}
\end{equation}
where $\Tmsc \in \{\Tmsc_{\rm A}, \Tmsc_{\rm L}\}$ is the tariff class.  Let $\mu_{\rm A}^*$  and $\mu_{\rm L}^*$  be the Ramsey pricing for the two-part tariff and linear tariff classes, respectively.

\subsection{Technology Adoption Model}
We now present a dynamical system  for  the adoption of distributed energy resources such as solar PV and energy storage.  We assume that the adoption decision of a residential customer is based on his investment's payback time, which  depends on the cost of solar PV and the reduced payment for consumption.
Instead of considering  individual adoption decisions, we model the adoption process for the entire customer population.

Following the standard innovation adoption theory \cite{Mahajan&Muller&Bass:91bookchap}, for a given tariff $T$ and per-unit (kWh) PV purchasing cost $\xi$, let the installed renewable capacity in aggregation be $s(t,T,\xi)$ at time $t$.  Illustrated in
Fig. \ref{fig_dif},  $s(t,T,\xi)$ is referred to as the {\em PV adoption curve} and  is defined by  the following equation:
\begin{equation}
s(t,T,\xi)=R_\infty(T,\xi) \eta(t),
\end{equation}
where $R_\infty(T,\xi)$ is the {\em market potential} of the PV adoption, and the cumulative installed fraction $\eta(t)$ is a sigmoid function satisfying $\eta(0)=0$ and $\mathop {\lim }\limits_{t \to \infty } \eta(t) =1$. The interpretation of the market potential is that, for fixed tariff $T$ and constant exogenous input $\xi$, the level of adoption eventually reaches $R_\infty(T,\xi)$. 

\begin{figure}[!h]
	\centering
	\begin{psfrags}
		\psfrag{t}[l]{$t$}
		\psfrag{s}[l]{$s(t,T,\xi)$}
		\psfrag{R}[r]{$R_\infty(T,\xi)$}
		\includegraphics[width=2in]{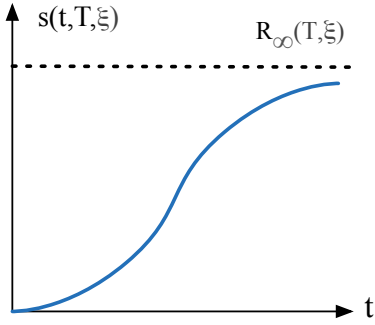}
		\vspace{-1em}
		\caption{Renewable adoption for fixed market potential.}
		\label{fig_dif}
	\end{psfrags}
\end{figure}
 
 Note that  the  shape of the adoption curve makes the model follow the so-called S-curve adoption of innovation.  This model has been used to model the adoption of renewable technology, and there is a parametric form of $R_\infty(T,\xi)$ that can be used in practice \cite{Beck:09}.  A well known form of $\eta(t)$ is from the Bass model \cite{Bass:69MS}.

%

%

The adoption curve  $s(t,T,\xi)$  does not capture the dynamics of the adoption process, however; it only describes the evolution of the adoption for {\em fixed }tariff $T$  and  PV cost $\xi$  throughout the adoption process.  In reality, the tariff is set by the utility periodically and the cost of PV declines. The evolution of the actual installed PV capacity in each period depends not only on the tariff and cost in that period but also on those in previous periods.  In other words, the installed PV capacity has to be calculated using not a single but a collection of such $S$-curves.  The dynamics of PV capacity evolution  is presented in Section~\ref{sec:III}.

\section{Dynamics and Stability of Adoption}\label{sec:III}
\subsection{Dynamics of Technology Adoption}
\label{sec3A}
We now introduce a discrete-time dynamical system model for the PV adoption process where the time index $k$ corresponds to the rate-setting epoch of the retailer.   The state $\sigma_k=(T_k, R_k)$ of the dynamical system includes the tariff $T_k$ set by the retailer at the beginning of the tariff period  and the installed PV capacity $R_k$ at the end of the tariff period.
The evolution of the system state is governed by the system equation
\begin{equation}
\sigma_{k+1}=f(\sigma_k,\chi_k),
\label{eq:dynamics}
\end{equation}
where $\chi_k=(\theta_k,\xi_k)$ is the exogenous (input) process containing the expected operating cost $\theta_k$ and the per-unit purchasing cost of PV $\xi_k$.  In analyzing the stability of the adoption process, we set the exogenous input to constant, $\chi_k=\chi$.  In general, the exogenous input can be time varying, especially when we  consider {\em controlled adoption} that sets tariff in response to varying costs. See case studies and qualitative results given in Section \ref{sec4B}.

The state evolution is assumed Markovian following
$$R_k \rightarrow T_{k+1}\rightarrow R_{k+1}$$ that corresponds to the decision process of the retailer who observes the level of adoption $R_k$ before setting the tariff $T_{k+1}$ for the next period.   Note that, at the beginning of period $k+1$, the installed PV capacity is $R_{k}$.  The installed capacity $R_{k+1}$ at the end of the period $k+1$ is obtained from the adoption curve associated with $T_{k+1}$  by $s(t_{k}+1,T_{k+1},\xi_{k})$ where $t_k$ is such that $s(t_k,T_{k+1},\xi_{k})=R_{k}$.  

A particularly relevant pricing policy is the myopic (greedy) Ramsey pricing that maximizes the consumer surplus subject to the break-even constraint in each period:
\begin{eqnarray}
\label{TfromR}
T_{k+1}&=& \mathop {\arg \max }\limits_{T\in \mathscr{T},\overline {\rm rs} (T,\theta_k,R_k) = 0} \overline {\rm cs} (T,R_k),\\
R_{k+1} &=& \left\{ \begin{array}{l}
{R_k},\quad \quad \quad \quad   {\rm{if   }}\;{R_\infty }({T_{k + 1}},\xi_{k} ) < {R_k};\\s(1 + {\eta ^{ - 1}}(\frac{{R_k}}{{R_\infty }({T_{k + 1}},\xi_{k} )}),{T_{k + 1}},\xi_{k} ),\;   {\rm o.w.}
\end{array} \right.
\label{Rdynamics}
\end{eqnarray}

Fig. \ref{fig_dyn} illustrates the evolution of the system state.  Suppose that at time $k$, the tariff is $T_k$ and the installed capacity is $R_k$ on the (blue) adoption curve $s(t,T_k,\xi_{k-1})$. Given $(R_k,T_k)$,  a new tariff $T_{k+1}$ is obtained for the next period, which has the (red) adoption curve $s(t,T_{k+1},\xi_k)$.  By shifting $s(t,T_{k+1},\xi_k)$ such that it intersects $s(t, T_k,\xi_{k-1})$ at $R_{k}$, we obtain the installed capacity $R_{k+1}$ for the $(k+1)$th period. 

\begin{figure}[!t]
	\centering
	\includegraphics[width=2.8in]{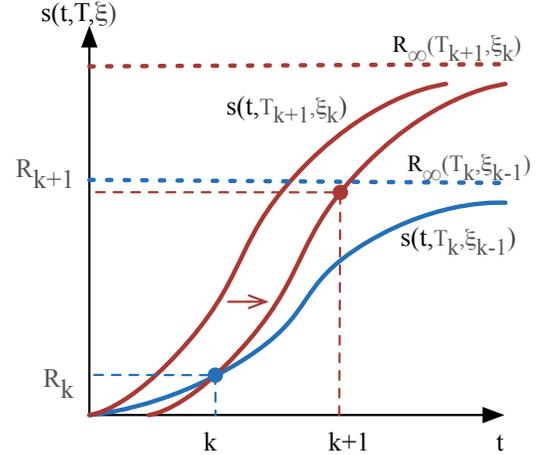}
	\vspace{-1em}
	\caption{Dynamics of renewable adoption when ${R_\infty }({T_{k + 1}},\xi_k )\geq {R_k}$.}
	\label{fig_dyn}
\end{figure}

A comment on the conditions in (\ref{Rdynamics}) is in order.  Typically, especially at the early stage of adoption, tariff $T_{k+1}$ will set the adoption curve to one with higher market potential as illustrated in Fig. \ref{fig_dyn}.  However, the situation that 
$R_\infty(T_{k+1},\xi_{k}) < R_k$ can happen when there is an exogenous shock such as a policy intervention such that $(T_{k+1},\xi_{k})$ leads to an adoption curve that has a lower market potential. For instance, the regulator may remove a certain tax incentive for PV adoption, which makes PV adoption less attractive and halts the adoption process.

\subsection{Death Spiral and its Existence Conditions}
The notion of death spiral is associated with the trajectory of a dynamical system defined through the tariff policy $\mu$ and the adoption curve.
\begin{Def}[Death spiral and critical adoption level]
	\label{def:cridif}
	An evolution of the dynamical system (\ref{TfromR}-\ref{Rdynamics}) states starting from $\sigma_0$ is a death spiral induced by tariff policy $\mu$ if it reaches a state $\sigma_{k_o}$ for which the optimization (\ref{TfromR}) to determine $T_{k_o+1}$ is not feasible.   The critical adoption level $R_\mu^\sharp$ is the supremum of $R$ at which a revenue adequate tariff exists
	\begin{equation}
	\label{eq:critical}
	R_\mu^\sharp=\sup\{R:\overline {\rm rs} (\mu(R,\theta),\theta,R)=0\}.
	\end{equation}
\end{Def}
We now establish a structural results on the critical adoption level. 

\begin{prop}
	\label{prop:rsharpdecrea}
	The critical adoption level $R^\sharp$ is monotonically decreasing on the retailer cost $\theta$. 
\end{prop}

We now focus on establishing existing conditions of death spiral. We assume exogenous parameters $\chi=(\theta,\xi)$  are fixed in this analysis.  For brevity, we drop notational dependencies on $\theta, \xi$, and $\chi$ when no confusion arises, and include them when the dependency on them plays a role in the analysis.

Our analysis relies on the  characterization of the {\em market potential function} defined as follows.
\begin{Def}[Market Potential function]
	Given a tariff policy $\mu$, the market potential function at adoption level $R$ is defined as
	\begin{equation}
	p_\mu(R)=R_\infty(\mu(R)).
	\end{equation}
\end{Def}
The market potential function serves as a surrogate for the more complicated iterative map $f$ in (\ref{eq:dynamics}). Being the maximum installation capacity on the adoption curve, $p_\mu (R)$ measures the headroom beyond the current installation capacity $R$.

The existence condition for death spiral is stated in Theorem \ref{th:deathcond} and illustrated in Fig. \ref{fig_deathdemo}.  Specifically, the death spiral occurs if the gap between $p_\mu(R)$ and $R$ is strictly positive in the left neighborhood of the critical adoption $R^\sharp$.

\begin{thm}[Existence condition of death spiral]	
	\label{th:deathcond}
	Given an initial state $\sigma_0$ with $R_0<R^\sharp$, a tariff policy $\mu$ generates a death spiral if there exists an $\epsilon>0$ such that
	\begin{itemize}
		\item $R_{k_0}\in(R^\sharp-\epsilon,R^\sharp]$ for some $k_0\geq 0$;
		\item $p_\mu(R)>R$ for all $R\in (R^\sharp-\epsilon,R^\sharp]$.
	\end{itemize}
	The condition is necessary and sufficient if $p_\mu(R)$ is monotonically increasing in $R$.
\end{thm}

\begin{figure}[!h]
	\centering
	\includegraphics[width=2in]{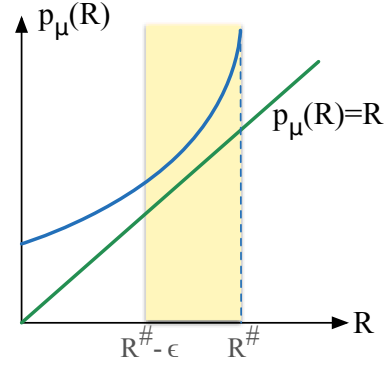}
	\vspace{-1em}
	\caption{Condition for death spiral as stated in Theorem \ref{th:deathcond}.}
	\label{fig_deathdemo}
\end{figure}

Theorem~\ref{th:deathcond} provides a way to check, at least numerically, the possibility of death spiral.  It is significant that the conditions in Theorem~\ref{th:deathcond} are necessary and sufficient when the market potential function is monotonically increasing as in most empirical cases studied in this paper.   The following Lemma provides a condition for the monotonicity of $p_\mu(R)$ with respect to $R$ and $\theta$.

\begin{lem} [Monotonicity of market potential function]
	\label{lem:monotonR}
The market potential function $p_\mu(R,\theta)$ for the Ramsey linear flat tariff $\mu^*_{\rm F}$ is monotonically increasing in $\theta$. It is also monotonically increasing in $R$ if, at the initial state when $R=0$,  the expected revenue from selling renewable in the retail market is higher than that in the wholesale market.

\end{lem}

The condition for the monotonicity of $p\mu(R,\theta)$ with respect to $R$ is natural since the retail price at the initial stage is almost always higher than the wholesale price, and the renewable incentive policy such as net-metering naturally ensures that the condition in Lemma \ref{lem:monotonR} holds.

We can now specialize Theorem \ref{th:deathcond} on the Ramsey linear tariff $\mu_{\rm L}^*$ that maximizes the consumer surplus subject to a revenue adequacy constraint.

\begin{corollary}[Death spiral condition for Ramsey linear tariff]
	\label{cor:deathram}
	Assume that for the Ramsey linear tariff $\mu^*_{\rm L}$, $p_{\mu_{\rm L}^*}(R,\theta)$ is monotonically increasing with respect to retailer cost $\theta$. There exists a threshold $\theta^\dagger$ such that $\mu_{\rm L}^*$ at 
	any retailer cost $\theta>\theta^\dagger$ induces a death spiral.
\end{corollary}

In particular, we can calculate the death spiral thresholds $\theta^\dagger$ for Ramsey linear flat tariff with some additive assumptions (see Corollary \ref{cor:deathflat} in the Appendix).

\subsection{Stable Adoption}
Death spiral is a form of instability.  We now consider conditions for {\em stable adoption}. In this context, we assume that both $\theta$ and $\xi$ are fixed.

We begin with standard definitions of the equilibrium and conditions for the stable equilibrium.
\begin{Def}[Stable equilibrium and stable adoption]\mbox{}For the state evolution defined in (\ref{eq:dynamics}), 
	\ben
	\item A state $\sigma^*$ is an equilibrium if $\sigma^*=f(\sigma^*)$.
	\item An equilibrium $\sigma^*$ is  Lyapunov stable if, for each $\epsilon>0$, there exists a $\delta=\delta(\epsilon)$ such that, for every trajectory $(\sigma_0,\sigma_1,\cdots)$ that is not a death spiral, $\left\| \sigma_0-\sigma^*\right\|<\delta$ implies  $\left\| \sigma_k-\sigma^*\right\|<\epsilon$ for all $k>0$.
	\item A trajectory $(\sigma_0,\sigma_1,\cdots)$ is a stable adoption if it converges to a stable equilibrium.
	\een
\end{Def}


\begin{thm}[Equilibrium, stability, and stable adoption]
	\label{th:stability}
	Given a tariff policy $\mu$ and a market potential function $p_\mu(R)$, if $R^*$ satisfies $p_\mu(R^*)\leq R^*$, then $\sigma^*=(\mu(R^*),R^*)$ is an equilibrium. Furthermore, 
	\ben
	\item  
	$\sigma^*$ is Lyapunov stable if there exists an $\epsilon>0$ such that $p_\mu(R)\leq R$ for all $R\in (R^*,R^*+\epsilon)$; 
	\item 
	A trajectory $(\sigma_0,\sigma_1,\cdots,\sigma_{k_0},\cdots)$ is a stable adoption with $\mathop {\lim }\limits_{k \to \infty } {\sigma _k} = {\sigma ^*}$ if we have $\forall R\in(R^*-\epsilon,R^*),$ $R<p_\mu(R)\leq R^*$   and $R_{k_0}\in (R^*-\epsilon,R^*)$ for some $k_0\geq 0$.
	\een 
\end{thm}

This equilibrium condition is intuitive; it states the case when the current level of installed PV capacity $R$ already reaches $R_\infty(\mu(R))$. If $p(R^*)<R^*$, from (\ref{Rdynamics}), $(\mu(R^*),R^*)$ is also an equilibrium. The first item in Theorem \ref{th:stability} provides the general existence condition of a stable equilibrium. Under this condition, a adoption reaching the right neighborhood of $R^*$ simply stays there. That said, for the case illustrated in Fig. \ref{fig_equidemo}, any capacity no less than the intersection is a stable equilibrium. In the second item, we specify a stable adoption that converges to $\sigma^*$. The graphical illustration is given in the area of Fig. \ref{fig_equidemo} shaded in yellow, with the intersection equal to the capacity of convergence.

\begin{figure}[!h]
	\centering
	\includegraphics[width=2in]{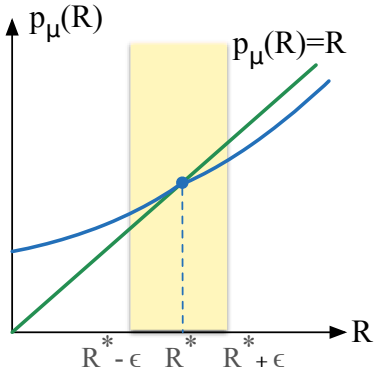}
	\vspace{-1em}
	\caption{Condition for stability as stated in Theorem \ref{th:stability}.}
	\label{fig_equidemo}
\end{figure}

\begin{thm}[Stable adoption via Ramsey two-part tariff]
	\label{th:ramseytwopart}
	For any initial state $\sigma_0=(T_0,R_0)$ with $0\leq R_0\leq p_{\mu^*_{\rm A}}(0)$, the Ramsey two-part tariff $\mu_{\rm A}^*$ induces a adoption approaching to the unique stable equilibrium $(\mu_{\rm A}^*(p_{\mu^*_{\rm A}}(0)),p_{\mu^*_{\rm A}}(0))$.
\end{thm}

Theorem~\ref{th:ramseytwopart} implies that the Ramsey two-part tariff can induce a stable adoption regardless of the values of $\theta$ and $\xi$.   It turns out that the solar adoption level at the stable equilibrium is quite low in such cases. The reason is that Ramsey two-part tariff $\mu_{\rm A}^*$ does not rely on markups in the volumetric price to achieve revenue adequacy, and it is the markups on volumetric prices that represent a strong incentive for PV adoption. Under the Ramsey pricing, the retail cost $\theta$ is recovered uniformly from the connection charges. Furthermore, the volumetric price of $\mu_{\rm A}^*$ does not depend on the solar integration capacity.

The following lemma provides a so-called critical connection charge that guarantees a stable adoption for two-part tariffs.

		\begin{lem}[Stable adoption via critical connection charge]
			\label{lem:criticalconnection}
			Let  $\theta^\sharp=\mathop {\min }\limits_R  (\mathop {\max }\limits_\pi (\mbbE[ (\pi-\lambda)^\top(D(\pi ,\omega ))-Rr_0(\omega)]))$. All two-part tariffs with  fixed connection charge $A^\sharp= (\theta-\theta^\sharp)/M$ induces a stable adoption.
			
		\end{lem}

Here, the max operation gives the maximum value of revenue a retailer extract from the consumers. From Definition \ref{def:cridif}, revenue at $R^{\sharp}$ matches the retailer's operating cost $\theta$. With the min operation, $\theta^\sharp$ means that, for any renewable adoption $R$, the retailer can always find a price $\pi$ such that its fixed cost $\theta^\sharp$ is covered. Therefore, by intuition, linear tariffs refrain from death spiral if the retailer cost is no more than $\theta^\sharp$.

\subsection{Limiting Diffusion Capacity}
In this subsection, we are interested in finding the highest level of PV adoption $R^\dagger$  achievable by a stable adoption.   The following definition formalizes the notion of limiting adoption capacity.

\begin{Def}[Limiting adoption capacity]   The limiting adoption capacity for a tariff class is the supremum of the installed capacity associated with a stable equilibrium achievable by a stable adoption.
\end{Def}

 Note that the highest consumer surplus can be achieved at the limiting adoption capacity (if the PV cost is considered as sunk cost and the tariff is properly designed). From (\ref{eq:critical}) we know that the critical adoption level for a tariff class is the same for all tariff policies (including the Ramsey tariff) that satisfy the break-even condition.   Let $R_{\mu_{\rm L}}^\sharp(\theta)$ be the critical adoption level for linear tariffs under retailer cost $\theta$. The following theorem characterizes the limiting adoption capacity for two-part tariffs satisfying the break-even constraint. Note that $R_{\mu_{\rm L}}^\sharp(\theta)$ is sufficient to include the critical adoption levels under two-part tariffs since (\ref{eq:rsexpression}) shows the connection charge $A$ and retailer cost $\theta$ are homogeneous in the expression of retail surplus.

\begin{thm}[Limiting capacity characterization]
	\label{deflim}
	
	Assume 1) $p(R_{\mu_{\rm L}}^\sharp(\theta))$ is increasing on $\theta$, and 2) $p(R_{\mu_{\rm L}}^\sharp(\theta))\geq p(R)$ for all $0\leq R \leq R_{\mu_{\rm L}}^\sharp(\theta)$ and $\theta\geq 0$. The limiting capacity $R^\circ$ for two-part tariffs subject to the break-even condition is equal to $R_{\mu_{\rm L}}^\sharp(\theta^\circ)$, where $\theta^\circ$ satisfies $p(R_{\mu_{\rm L}}^\sharp(\theta^\circ))=R_{\mu_{\rm L}}^\sharp(\theta^\circ)$.
\end{thm}

\begin{figure}[!h]
	\centering
	\includegraphics[width=2.5in]{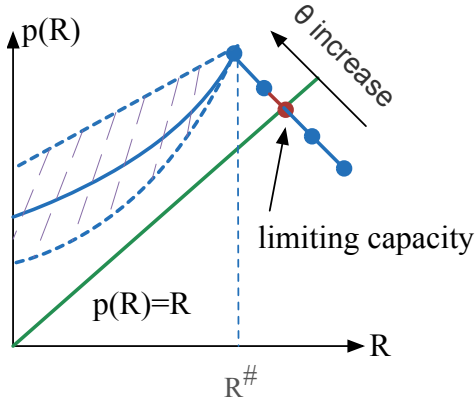}
	\vspace{-1em}
	\caption{Limiting capacity for two-part tariffs as stated in Theorem \ref{deflim}.}
	\label{fig_critical}
\end{figure}

Fig. \ref{fig_critical} gives an illustration for Theorem \ref{deflim}. The shaded part represents all break-even linear tariff policies under a particular retailer cost. The blue dots represent the critical adoption levels when the retailer cost changes. The red dot at the intersection is the limiting capacity. To achieve the  limiting capacity requires a complicated control of the tariff policy. For example, the tariff policy may require constant changes of the connection charge and the volumetric price which are impractical in practice. 

Currently we cannot characterize a specific tariff policy that achieves the limiting capacity. However, one can obtain a \textit{capacity lower bound} for the limiting capacity by restricting tariff policies to those with a fixed connection charge and volumetric price $\mu_{\rm L}^*$. Assume that the market potential function $p(R,\theta)$ of the Ramsey linear tariff $\mu^*_{\rm L}$ is monotonically increasing with respect to retailer cost $\theta$. If  $\mu^*_{\rm L}$ induces a death spiral,   by adding the minimum connection charge (limiting connection charge) so that there is a stable adoption, we can find the maximum adoption under such tariff policies at the equilibrium (We name this adoption level as the capacity lower bound).

 Furthermore, if the volumetric price is flat ($\mu_{\rm F}^*$), and the assumption in Corollary \ref{cor:deathflat} holds,  the capacity lower bound is $R^\dagger$ and the fixed connection charge that achieves the capacity lower bound is $A^\dagger=(\theta-\theta^\dagger)/M$.

\section{An Empirical Study}

We analyzed renewable adoption dynamics in both short-run and long-run cases within a hypothetical distribution utility facing
the wholesale price  and residential demand in New York City.  The same settings of linear demand model, consumption profile, revenue estimation, and solar PV data weer used as in \cite{Munoz-Alvarez&Tong:17TPS1,Munoz-Alvarez&Tong:17TPS2}.

The default tariff of the Consolidated Edison Company of New York (ConEd) in 2015 for its 2.2 million residential customers was a two-part tariff $T^{\rm CE}$ with a connection charge $A^{\rm CE}$ equal to $\$0.52/{\rm day}$ and a flat volumetric price $\pi^{\rm CE}$ equal to $\$0.172/{\rm kWh}$. We substituted the prevailing retail tariff $T^{\rm CE}$ into the utility's break-even condition in (\ref{eq:ramplan}) to estimate the utility's daily fixed costs, which amount to $\theta^{\rm CE}=\$6.03 {\rm M}$. A consumer surplus of $\overline{\rm cs}_0(T^{\rm CE})=\$9.54{\rm M}$ was assumed as in \cite{Munoz-Alvarez&Tong:17TPS1}.

The integration of solar PV was modeled based on a simulated 5kW-DC-capacity rooftop system in NYC. The market potential $R_\infty$ was computed based on the expected payback years $t^{\rm PB}=\xi/\mathbb{E}[\pi^\top r_0(\omega)]$ at the time of purchasing. We took the solar PV cost of NYC in 2015 as the initial solar cost $\xi_0=\$4250/{\rm kW}$\footnote{The solar cost data in New York State starting from 2009 can be found at https://www.nysolarmap.com/}. An exponential fit in \cite{Beck:09} was adopted in calculating the market potential: $R_\infty=R^{\rm MS}\cdot e^{-0.3t^{\rm PB}}$. As in \cite{NREL:09}, the total market size $R^{\rm MS}$ was set to be $90\%$ of all customers installing, and $\eta(t)$ was set to model a medium-rate adoption using the Bass model.

\subsection{Short-run Analysis}

In the short-run analysis, exogenous parameters including the retailer's cost and the solar cost were fixed during the evolution of states. We examined both flat and dynamic tariffs.

\begin{figure}[!t]
	\centering
	\includegraphics[width=3in]{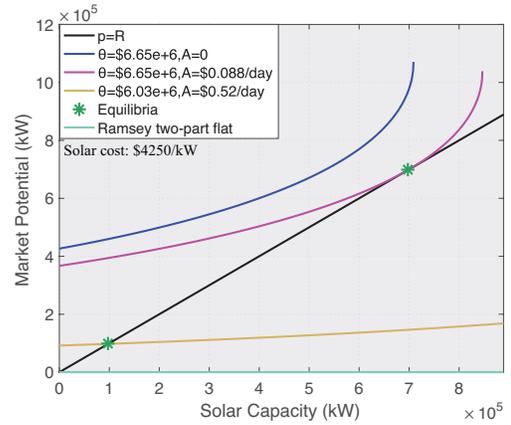}
	\vspace{-1em}
	\caption{Potential function in short run analysis.}
	\label{fig_srp2r}
\end{figure}

Fig. \ref{fig_srp2r} shows the market potential functions under different  tariff policies with flat volumetric prices. For each tariff class, the market potential function was increasing with respect to solar capacity (The market potential function of Ramsey two-part tariff was horizontal). The equilibrium adoption capacity of the Ramsey two-part tariff $\mu^*_{\rm A}$ was almost at 0, which stalled the solar adoption (green curve). This phenomenon was due to the low volumetric price under such tariff policy, leading to long payback time. The two-part tariff with volumetric price $\mu_{\rm F}^*$ and fixed connection charge $A^{\rm CE}$ (as currently used by ConEd), induced a stable equilibrium with the solar capacity equal to 97.7MW (brown curve).

When the retailer cost was increased by about 10\% to \$6.65M, Ramsey linear flat tariff $\mu^*_{\rm F}$ induced a death spiral (blue curve).   When a connection charge $A\geq \$0.088/{\rm day}$ was introduced, the tariff stayed off the death spiral and achieved a stable equilibrium.   Moreover, when the connection charge $A^\dagger= \$0.088/{\rm day}$ (less than 20\% of the current ConEd's) was imposed, the capacity lower bound $R^\dagger$ could be achieved (magenta curve) at 0.7M kW (the limiting capacity is 1M kW). A connection charge higher than $A^\dagger$ will induce a stable adoption. The difference between $A^\dagger$ and the current connection charge imposed by ConEd revealed that the current connection charge not only avoided the death spiral, but also left an adequate stability margin even when the retailer's operating cost rose by 10\%.

\begin{figure}[!t]
	\centering
	\includegraphics[width=3in]{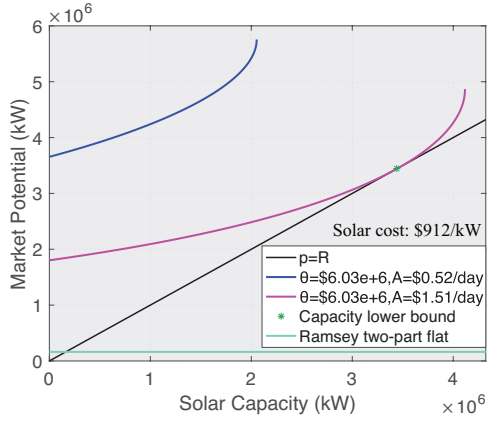}
	\vspace{-1em}
	\caption{Potential functions of Ramsey two-part tariff $\mu^*_{\rm A}$, linear tariff $\mu^*_{\rm F}$, and two-part tariff with volumetric price $\mu_{\rm F}^*$ and limiting connection charge. }
	\label{fig_srp2r2}
\end{figure}

Similar market potential function curves with a lower solar cost are shown in Fig. \ref{fig_srp2r2}. Lowering solar cost by more than 75\% leads to a  significantly higher level of adoption.  The current ConEd tariff at \$0.51/day that induced stable adoption when the solar cost was high induced death spiral.  The death spiral can be mitigated by increasing the connection charge to \$1.51/day.

Fig. \ref{fig_p2rdyn} shows the market potential functions under different tariff policies with time-varying volumetric prices. The market potential function might not be increasing on solar capacity in this case (while the market potential function of Ramsey two-part tariff is still horizontal). The adoption equilibrium of the Ramsey two-part tariff $\mu^*_{\rm A}$ was almost at 0, which stalled the solar adoption (green curves) as well. The critical connection charge ($A^\sharp=\$0.97/{\rm day}$), regardless of solar costs, induced stable adoptions (brown curves) which have higher equilibrium capacity than under $\mu^*_{\rm A}$. 

\begin{figure}[!t]
	\centering
	\includegraphics[width=3in]{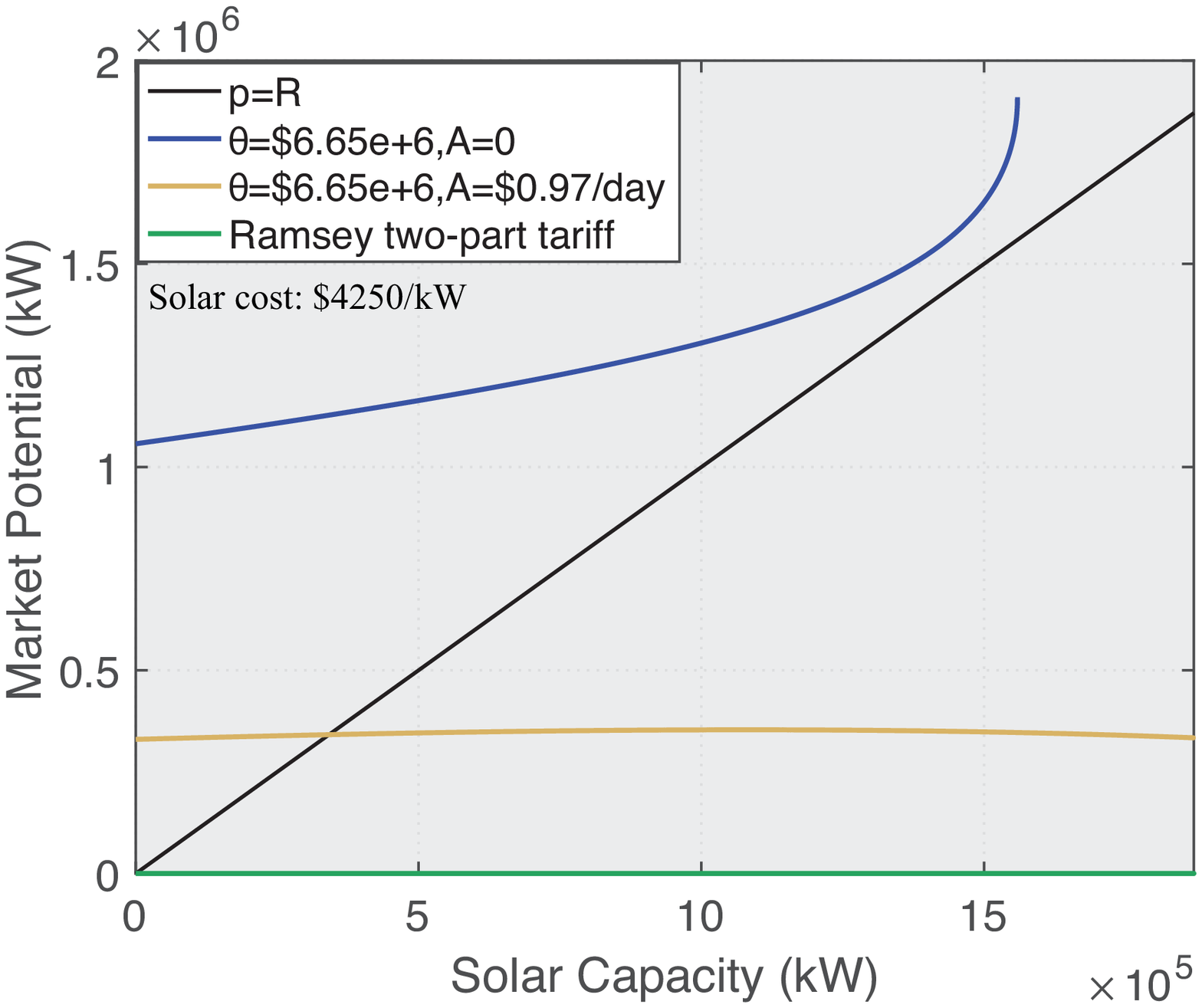}
	\includegraphics[width=3in]{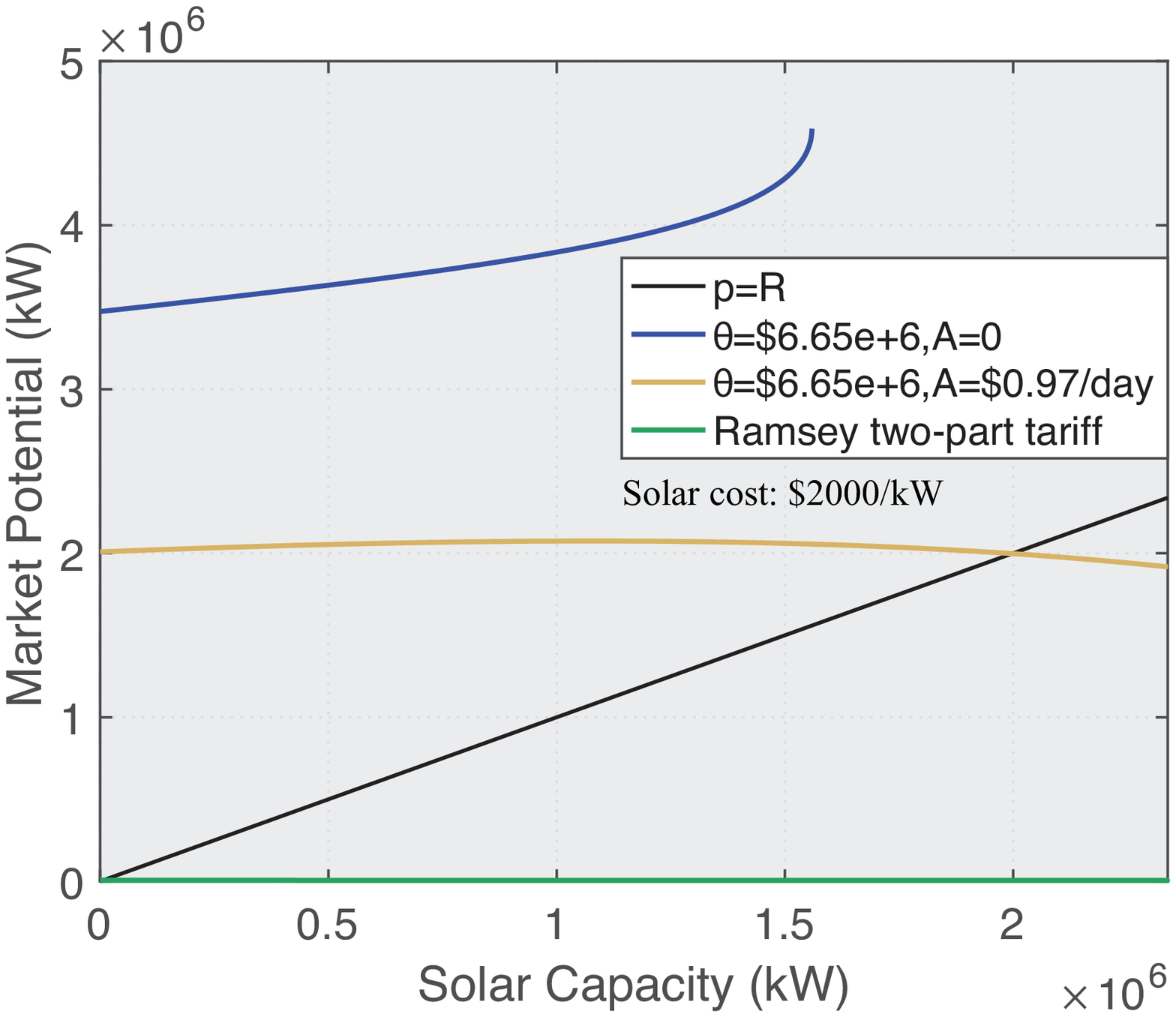}
	\vspace{-1em}
	\caption{Potential functions of Ramsey two-part tariff $\mu^*_{\rm A}$, linear tariff $\mu^*_{\rm L}$, and two-part tariff with volumetric price $\mu_{\rm L}^*$ and critical connection charge.}
	\label{fig_p2rdyn}
\end{figure}


\subsection{Long-run Analysis}
\label{sec4B}

Fig. \ref{fig_srsw2t} shows a comparison of total consumer surplus between two two-part tariffs over 20 years. One tariff policy was the Ramsey two-part tariff $\mu^*_{\rm A}$, the other was the two-part tariff with the connection charge $A^\dagger$=\$1.51/day and volumetric price $\mu_{\rm F}^*$.  In both cases, the solar cost (\$912/kW) and retailer cost were fixed.

\begin{figure}[!t]
	\centering
	\includegraphics[width=3in]{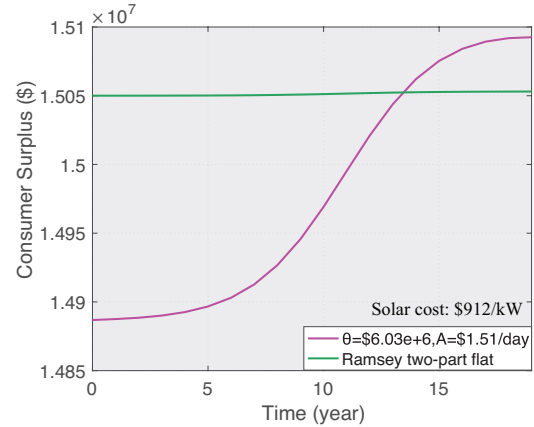}
	\vspace{-1em}
	\caption{Consumer surplus induced by Ramsey two-part tariff $\mu^*_{\rm A}$ and two-part tariff with limiting connection charge.}
	\label{fig_srsw2t}
\end{figure} 

As shown in Fig. \ref{fig_srsw2t}, the total consumer surplus induced by $\mu^*_{\rm A}$  had slow growth. Under the two-part tariff with $A^\dagger$, the consumer surplus was low at first but eventually became higher than under $\mu^*_{\rm A}$ due to a higher solar installation. This comparison illustrated the trade-off in achieving long-run consumer surplus optimization: the need to add connection charge for cost recovery and the need to limit the connection charge to promote PV adoption. The Ramsey two-part tariff $\mu^*_{\rm A}$, which maximizes the consumer surplus greedily, is not the optimal choice for consumer surplus maximization in the long run.

Fig. \ref{fig_bar3} illustrates effects of connection charge on the PV adoption capacity and the consumer surplus. As the connection charge increased from the minimum (limiting) connection charge of \$1.51/day to the maximum (Ramsey) connection charge of \$2.74/day, the adoption capacity decreased for all years.   The lower panel of Fig. \ref{fig_bar3} shows that the limiting connection charge does not necessarily lead to the highest consumer surplus, which was achieved at the fixed connection charge of $A$=\$1.78/day. 

 In this context, total consumer surplus does not fully capture the effect of solar cost.  We define, for the long run numerical studies, a (narrow) notion of social welfare that includes the cost of PV, specifically, 
		$\overline {\rm sw}=\sum_{year} \overline {\rm cs}-R\cdot\xi$.
We show how  social welfare changes over the years in Fig. \ref{fig_swcum}. Although Fig. \ref{fig_bar3} has shown the increase of consumer surplus caused by more PV capacity, Fig. \ref{fig_swcum} reveals that the highest social welfare achieves at $A$=\$2.74/day, which means the increase of aggregated surplus brought by PV installation cannot cover the cost of solar within the time scope of this case (20 years). However, this situation may change if we look at a longer time scale, due to the surplus superiority and saturation of new PV installation.

\begin{figure}[!t]
	\centering
	\includegraphics[width=3in]{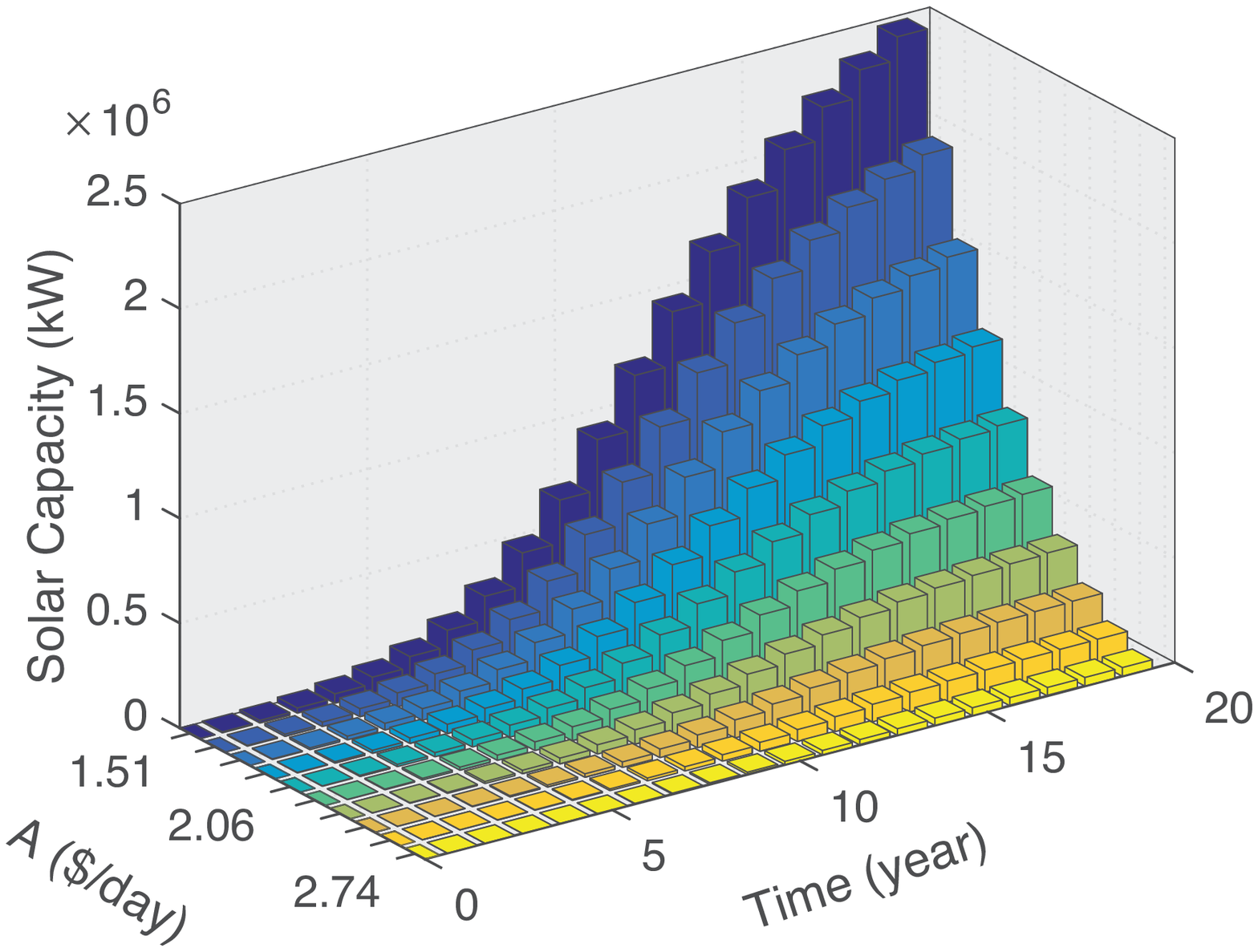}
	\includegraphics[width=3in]{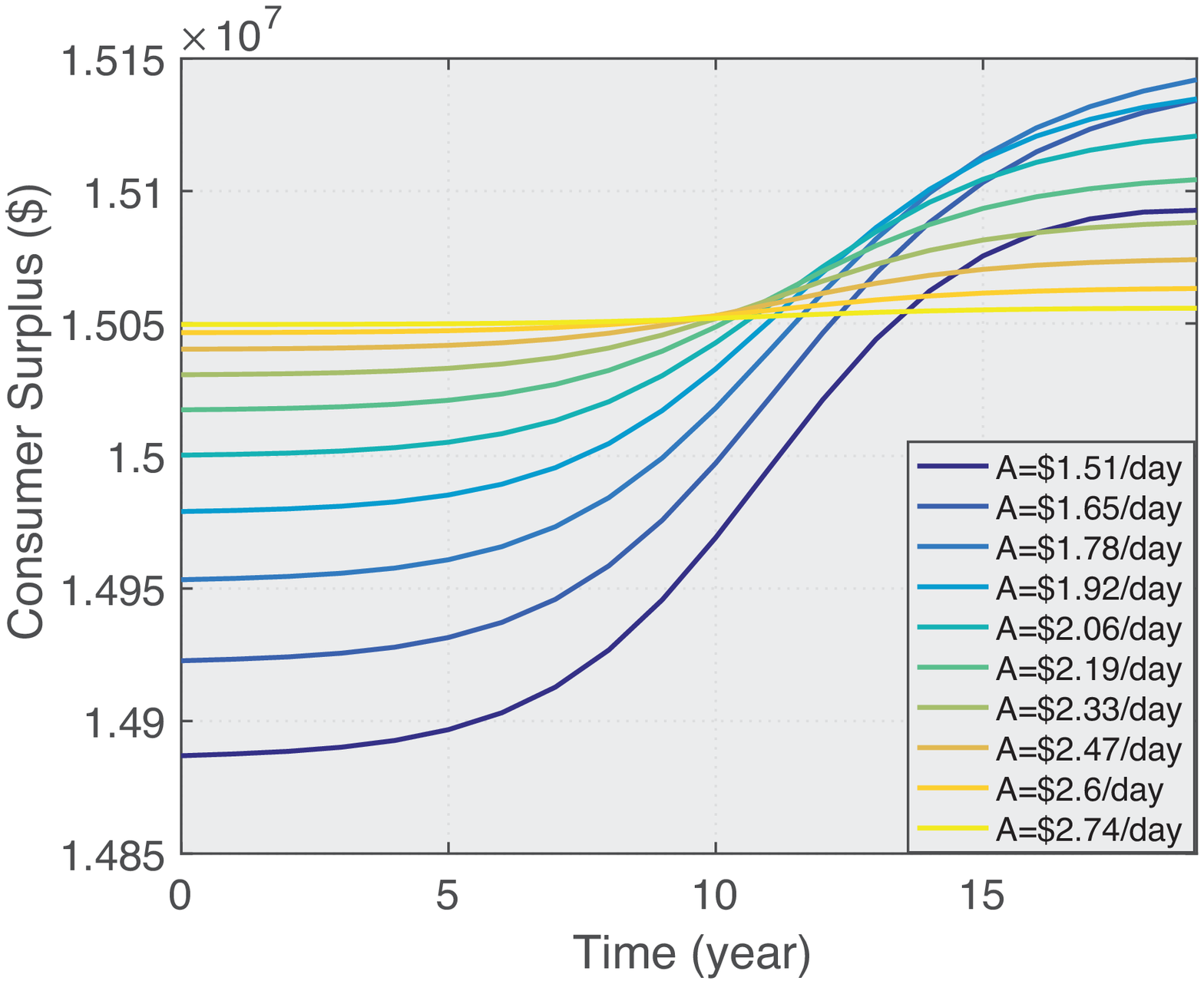}
	\vspace{-1em}
	\caption{ Solar capacity and consumer surplus of stable PV adoptions under different connection charges.}
	\label{fig_bar3}
\end{figure}

\begin{figure}[!h]
	\centering
	\includegraphics[width=3in]{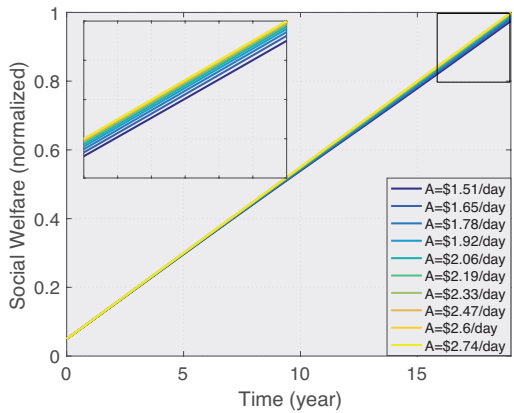}
	\vspace{-1em}
	\caption{Cumulative social welfare under different connection charges.}
	\label{fig_swcum}
\end{figure}

Fig. \ref{fig_r2tvtheta} and Fig. \ref{fig_r2tvzeta} show the long-run solar adoption dynamics (flat volumetric prices) under processes of increasing retailer cost and decreasing solar cost. In both cases, death spiral was induced under the same tariff. Adding limiting connection charges, however, could stay off the death spiral and achieve stable adoptions (magenta curve). Moreover, while introducing limiting connection charges lowered the speed of solar integration, its adoption capacity was higher than under fixed connection charge in the long run. The fixed connection charge case generates a death spiral which stalls further solar adoption.

Although it is beyond the scope of this paper to analyze the adoption process with varying costs, some insights could be gained from the numerical studies.  Fig. \ref{fig_r2tvtheta} shows that increasing retailer costs does not affect the capacity lower bound (dashed magenta line). If the connection charge keeps unchanged (blue curve), the death spiral is induced, even faster than under fixed retailer cost $\theta$. The reason is that increasing $\theta$ not only accelerates PV adoption by lifting market potential but also lowers the critical adoption level $R^\sharp$. In Fig. \ref{fig_r2tvzeta}, we show that decreasing solar costs does not change the critical adoption level (blue dashed line). The death spiral still comes earlier because of the lifting market potential. The capacity lower bound rises with solar cost decrease. The effect of increasing wholesale prices (not plotted here) is similar to a combination of the above. The death spiral occurs faster because of the growing market potential and lowered critical adoption level. The capacity lower bound also increases in this case.

\begin{figure}[!t]
	
	\centering
	\includegraphics[width=3in]{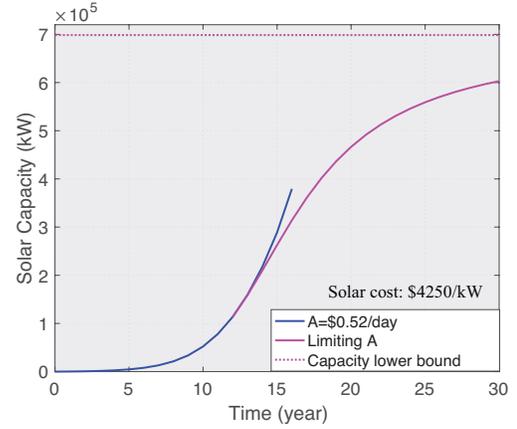}
	\vspace{-1em}
	\caption{Long-run solar adoption with retailer cost increasing by 2\% every year from $\theta^{\rm CE}$.}
	\label{fig_r2tvtheta}
\end{figure}

\begin{figure}[!t]
	\centering
	\includegraphics[width=3in]{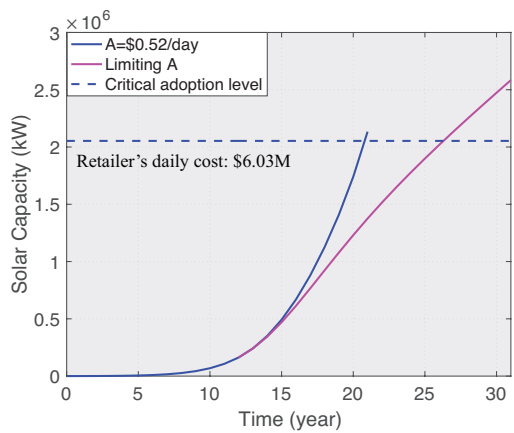}
	\vspace{-1em}
	\caption{Long-run solar adoption with solar cost decreasing by 5\% every year from $\xi_0$.}
	\label{fig_r2tvzeta}
\end{figure}

\section{Conclusion}\label{sec:VI}

The goal of developing an analytical framework for characterizing the dynamics of PV adoption is to obtain insights into the behavior of the adoption process and roles of tariff on adoption.  Based on the analysis, several conclusions can be made.  First, linear tariffs are prone to death spiral when the fixed cost of the utility rises beyond a certain level.  Our model also allows us to quantify situations when death spiral is imminent.  Second,  a small connection charge has the effect of stimulating adoption whereas a high connection charge tends to inhibit adoption.  
Our analysis provides a way to set the right level of connection charge to induce a stable adoption process and greater adoption capacity.

Finally, we comment on some of the limitations of this paper and  future work. The current analytical model does not include so-called partial net-metering as proposed by some public utility commissions. The difference of such a mechanism with net-metering is that it compensates the electricity sent back to the grid at a price lower than the retail price. The changed pricing structure thus significantly increases the difficulty of integrating the consumer decision model into adoption dynamics. It might be interesting to look at how PV adoption changes when the compensation price varies between the retail price and the wholesale price.

Another issue we have not addressed is the design of retail tariff that achieves both incentive compatibility and fairness in PV adoption. The current practice of net-metering boosts the surplus of PV owners at the expense of non-PV owners. Possible solutions may include differentiated pricing for consumers with/without PV installations and (load profile) data-driven tariff design.

\section*{Acknowledgment}

The authors would like to thank Steven H. Low from Caltech and Desmond Cai from Institute of High Performance Computing for insightful discussions.

{
\bibliographystyle{ieeetran}
\bibliography{BIB}
}

\section*{Appendix}
\begin{prop}
	\label{prop:rsmax}
	For a linear tariff $\pi^\sharp=\mu(R^\sharp,\theta)$ at the critical adoption capacity $R^\sharp$, we have $$\pi^\sharp=\arg \mathop {\max }\limits_\pi \overline {\rm rs} (T,\theta,R^\sharp).$$ Denote $$\overline {\rm rs}^{\rm M}(A,\theta,R^\sharp)=\mathop {\max }\limits_\pi \overline {\rm rs} (T,\theta,R^\sharp).$$ The relation among $T,\pi$, and $A$ are defined in Section II-A.
\end{prop}
\begin{proof}[Proof:]
	If there exists a $\pi_0$ such that $$\overline {\rm rs} (T_0,\theta,R^\sharp)>\overline {\rm rs} (T^\sharp,\theta,R^\sharp)=0,$$ there must exist $R^\prime>R^\sharp$ satisfying $\overline {\rm rs} (T_0,\theta,R^\prime)>0$ due to the continuity. Thus a contradiction is induced with the definition of critical adoption level.
\end{proof}

\begin{proof}[Proof of Proposition \ref{prop:rsharpdecrea}:]
	Leveraging Proposition \ref{prop:rsmax}, for a retailer cost $\theta_1$ and the corresponding critical adoption level $R^\sharp_1$, we have $\overline {\rm rs}^{\rm M}(A,\theta_1,R^\sharp_1)=0$. Hence, with the expression of $\overline {\rm rs}$ in (\ref{eq:rsexpression}), we have $\overline {\rm rs}^{\rm M}(A,\theta_2,R^\sharp_1)<0$ for all $\theta_2>\theta_1$. Thus $R^\sharp_2<R^\sharp_1$.
\end{proof}

\begin{prop}
	\label{prop:increase}
	For a trajectory $(\sigma_0,\sigma_1,\cdots)$, if $p(R_k)>R_k$, we have $R_k<R_{k+1}=h(R_k,\mu(R_k,\theta))<p(R_k)$. 
\end{prop}
\begin{proof}[Proof:]
It directly holds from Equation (\ref{Rdynamics}).
\end{proof}

\begin{prop}
	\label{prop:converge}
	If there is an $\epsilon>0$ such that $R<p(R)\leq R^*$ for all $R\in(R^*-\epsilon,R^*)$ with $p(R^*)=R^*$, then for each $R_0\in (R^*-\epsilon,R^*)$, we have $\mathop {\lim }\limits_{k \to \infty } {R _k} = {R ^*}$.
\end{prop}
\begin{proof}[Proof:]
	Leveraging Proposition \ref{prop:increase}, $\{R_t\}$ is strictly increasing and bounded by $R^*$. We suppose that $\{R_t\}$ converges to $R^\prime\in (R_0,R^*)$. It can be induced that $h(R^\prime,\mu(R^\prime,\theta))=R^\prime$. As $p(R^\prime)>R^\prime$, there is a contradiction with Proposition \ref{prop:increase}. Hence $\{R_k\}$ must converge to $R^*$ (Monotone convergence theorem).
\end{proof}

\begin{proof}[Proof of Theorem \ref{th:deathcond}]
	Sufficiency $\Rightarrow$: Leveraging Proposition \ref{prop:increase}, $\{R_k\}$ is monotonically increasing. Suppose $R^\sharp$ is an upper bound of $\{R_k\}$.  Thus there exists an $R^\prime\in (R_{k0},R^\sharp]$ such that $\{R_k\}$ converges to $R^\prime$ (Monotone convergence theorem). Hence, $h(R^\prime,\mu(R^\prime,\theta))=R^\prime$. As $p(R^\prime)>R^\prime$, there is a contradiction with Proposition \ref{prop:increase}. Thus $R^\sharp$ is not an upper bound of $\{R_k\}$, indicating that the death spiral occurs.
	
	If $p(R)$ is monotonically increasing, the necessity can also be proved.
	
	Necessity $\Leftarrow$: Since a death spiral is induced, there must exist $R_0\leq R_{k1}<R^\sharp$ such that $p(R_{k1})>R^\sharp$ (Otherwise $R_{k+1}<p(R_k)\leq R^\sharp$ for all $k$, indicating there is no death spiral). Moreover, as $p(R)$ is monotonically increasing, $p(R)>p(R_{k1})>R^\sharp>R$ holds for $R\in(R_{k1},R^\sharp]$. Thus the necessity is proved.
\end{proof}

\begin{proof}[Proof of Lemma \ref{lem:monotonR}]
	The volumetric price of the Ramsey linear flat tariff is characterized by
	\begin{equation}
	\pi _{\rm{F}}^* = \frac{{\mathbb{E}[{\lambda ^ \top }\partial D({\bf{1}}\pi _{\rm{F}}^*,\omega )/\partial \pi _{\rm{F}}^*]}}{{{{\bf{1}}^ \top }\mathbb{E}[\partial D({\bf{1}}\pi _{\rm{F}}^*,\omega )/\partial \pi _{\rm{F}}^*]}} - \frac{{\gamma  - 1}}{\gamma }\frac{{{{\bf{1}}^ \top }\mathbb{E}[D - r(\omega )]}}{{{{\bf{1}}^ \top }\mathbb{E}[\partial D/\partial \pi _{\rm{F}}^*]}}
	\label{eq:ramseyflatenergyprice}
	\end{equation}
	where $\gamma$ is the Lagrange multiplier of (\ref{eq:ramplan}), with $\frac{\gamma-1}{\gamma}\in[0,1]$. The revenue adequacy constraint can be reformulated as
	\begin{equation}
	(\mathbf{1}\pi_{\rm F}^*)^\top \mathbb{E}[D-r(\omega)]=\theta+\mathbb{E}[\lambda^\top (D-r(\omega))]
	\label{eq:revadeflat}
	\end{equation}
	Differentiating both sides of (\ref{eq:revadeflat}) over $\theta$	yields (after some deduction)
	\begin{equation}
	\frac{\partial \pi _{\rm{F}}^*}{\partial \theta}\Bigg(\mathbb{E}[{({\bf{1}}\pi _{\rm{F}}^* - \lambda )^ \top }\frac{{\partial D}}{{\partial \pi _{\rm{F}}^*}}] +\mathbf{1}^\top \mathbb{E}[ D - R{r_0}(\omega ) ]  \Bigg)=1
	\label{eq:partialtheta}
	\end{equation}
	Substituting (\ref{eq:ramseyflatenergyprice}) into (\ref{eq:partialtheta}) yields
	\begin{equation}
	\frac{1}{\gamma} \mathbf{1}^\top \mathbb{E}[D - R{r_0}(\omega )]  \frac{\partial \pi _{\rm{F}}^*}{\partial \theta} =1
	\end{equation}
	It is then clear $\frac{\partial \pi _{\rm{F}}^*}{\partial \theta}>0$, thus $p(R,\theta)$'s monotonicity in $\theta$ is proved.
	
	Differentiating both sides of (\ref{eq:revadeflat}) over $R$ yields (after some deduction)
	\begin{equation}
	\begin{array}{l}
	\mathbb{E}[{({\bf{1}}\pi _{\rm{F}}^* - \lambda )^ \top }\frac{{\partial D}}{{\partial \pi _{\rm{F}}^*}}\frac{{\partial \pi _{\rm{F}}^*}}{{\partial R}}] + {{\bf{1}}^ \top }\mathbb{E}[D - R{r_0}(\omega )]\frac{{\partial \pi _{\rm{F}}^*}}{{\partial R}}\\
	\quad \quad \quad \quad \quad \quad \quad \quad \quad \; = \mathbb{E}[{({\bf{1}}\pi _{\rm{F}}^* - \lambda )^ \top }{r_0}(\omega )]
	\end{array}
	\label{eq:partialR}
	\end{equation}
	Substituting (\ref{eq:ramseyflatenergyprice}) into (\ref{eq:partialR}) yields
	\begin{equation}
	\frac{1}{\gamma} {{\bf{1}}^ \top }\mathbb{E}[D - R{r_0}(\omega )]\frac{{\partial \pi _{\rm{F}}^*}}{{\partial R}}= \mathbb{E}[{({\bf{1}}\pi _{\rm{F}}^* - \lambda )^ \top }{r_0}(\omega )]
	\end{equation}
	Leveraging $\mathbb{E}[{({\bf{1}}\pi _{\rm{F}}^*|_{R=0} - \lambda )^ \top }{r_0}(\omega )]\geq 0$, it is then clear  $\frac{{\partial \pi _{\rm{F}}^*}}{{\partial R}}\geq 0$, thus $p(R,\theta)$'s monotonicity in $R$ is proved.
	
\end{proof}

\begin{proof}[Proof of Corollary \ref{cor:deathram}]
	
	It is not hard to observe that as $\theta\to +\infty$, the critical adoption level $R^\sharp \to 0$. Hence, with Proposition \ref{prop:rsharpdecrea}, there must exist $\theta^\dagger$ such that the market potential function $p(R,\theta)$ meets the condition in Theorem \ref{th:deathcond}. As we assume $p(R,\theta)$ is monotonically increasing with respect to retailer cost $\theta$ for the Ramsey linear tariff $\mu^*_{\rm L}$, a retailer cost $\theta>\theta^\dagger$ will still induce a death spiral.

\end{proof}

\begin{prop}
	\label{prop:piconvex}
	For the Ramsey linear flat tariff, if (i) the market potential $R_\infty(\cdot)$ is convex on the flat volumetric price $\pi_{\rm  F}^*$, (ii) consumers' demand function is affine with negative slope and random disturbance, \ie $D(\pi,\omega)=B(\omega)-G\pi$, where $B(\omega)$ is the additive disturbance and $G$ positive semidefinite, and (iii) the assumption in Lemma \ref{lem:monotonR} holds, the market potential function $p(R,\theta)$ is not only monotonically increasing but also convex in $\theta$ and in $R$. 
\end{prop}

\begin{proof}[Proof:]
	Solving (\ref{eq:ramplan}) yields
	\begin{equation}
	\pi_{\rm F}^*(R,\theta)=\frac{-b(R)-\sqrt{b(R)^2-4ac(R,\theta)}}{2a}
	\end{equation}
	where $a= \mathbf{1}^\top G\mathbf{1}$,
	$b(R)=-\bar{\lambda}^\top G\mathbf{1}-\mathbf{1}^\top(\mathbb{E}[B(\omega)]-R\bar{r}_0)$, and $c(R,\theta)=\theta+\mathbb{E}[\lambda^\top(B(\omega)-R r_0(\omega))]$.
	
	Since the monotonicity directly holds from Lemma \ref{lem:monotonR}, we only need to show the convexity. As we have assumed $R_\infty(\cdot)$ to be convex on the the flat volumetric price $\pi_{\rm  F}^*$, we then need to prove $\pi_{\rm F}^*(R,\theta)$'s convexity on $\theta$ and $R$.
	
	\textbf{a)} On $\theta$: Differentiating twice $\pi_{\rm  F}^*(R,\theta)$ with respect to $\theta$ we have
	\begin{equation}
	\begin{array}{l}
	\frac{{d\pi_{\rm  F}^*}}{{d\theta }} = \frac{{ - 1}}{{2a}}( - 4a) \cdot \frac{1}{{\sqrt {b{{(R)}^2} - 4ac(R,\theta )} }}\\
	= \frac{2}{{\sqrt {{b^2} - 4ac} }} > 0
	\end{array}
	\end{equation}
	\begin{equation}
	\frac{d\pi_{\rm  F}^*}{d\theta^2}=\frac{2a}{\sqrt{b^2-4ac}(b^2-4ac)}
	\end{equation}
	Since $a=\mathbf{1}^\top G\mathbf{1}$ and $G$ positive definite, $\frac{d\pi^*}{d\theta^2}\geq 0$. 
	
	\textbf{b)} On $R$: differentiating $\pi_{\rm  F}^*(R,\theta)$ with respect to $R$ we have
	\begin{equation}
	\begin{array}{*{20}{l}}
	\begin{array}{l}
	{\pi_{\rm  F}^*(R) }^\prime = \frac{1}{{2a}}[ - {b^\prime } - (b{b^\prime } - 2a{c^\prime }) \cdot \frac{1}{{\sqrt {{b^2} - 4ac} }}]\\
	= \frac{1}{{2a\sqrt {{b^2} - 4ac} }}[2a{c^\prime } - b{b^\prime } - {b^\prime }\sqrt {{b^2} - 4ac} ]
	\end{array}\\
	{ = \frac{1}{{\sqrt {{b^2} - 4ac} }}[\frac{{( - b - \sqrt {{b^2} - 4ac} )}}{{2a}}{b^\prime } + {c^\prime }] = \frac{1}{{\sqrt {{b^2} - 4ac} }}(\pi_{\rm  F}^*(R){b^\prime } + {c^\prime })}
	\end{array}
	\end{equation}
	where: $b^\prime=\mathbf{1}^\top\bar{r}_0$ and $c^\prime=-\mathbb{E}[\lambda^\top r_0(\omega)]$.
	
	We differentiate twice $\pi_{\rm  F}^*(R)$
	
	\begin{equation}
	{\pi_{\rm  F}^*(R) }^{\prime\prime}=\frac{1}{2a}\frac{(b^\prime b-2ac^\prime)^2-(b^2-4ac)b^{\prime 2}}{\sqrt{b^2-4ac} (b^2-4ac)}
	\end{equation}
	With ${\pi_{\rm  F}^*(R) }^\prime\geq 0$, we have $2a{c^\prime } - b{b^\prime } - {b^\prime }\sqrt {{b^2} - 4ac} \geq 0$, which yields $(b^\prime b-2ac^\prime)^2-(b^2-4ac)b^{\prime 2}\geq 0$. Thus ${\pi_{\rm  F}^*(R) }^{\prime\prime} \geq0$ holds. 
	
\end{proof}

\begin{corollary}[Death spiral threshold for Ramsey linear flat tariff]	
	\label{cor:deathflat}
	For the Ramsey linear flat tariff $\mu^*_{\rm F}$, if (i) the market potential $R_\infty(\cdot)$ is convex on the flat volumetric price $\pi_{\rm  F}^*$, (ii) consumers' demand function is affine with negative slope and random disturbance, \ie $D(\pi,\omega)=B(\omega)-G\pi$, where $B(\omega)$ is the additive disturbance and $G$ positive semidefinite, and (iii) the assumption in Lemma \ref{lem:monotonR} holds, then the minimum of such thresholds as specified in Corollary \ref{cor:deathram} gives
	\begin{equation}
	\begin{array}{l}
	{\theta ^\dag } = \frac{1}{{4\mathbf{1}^\top G\mathbf{1}}}\Bigg[b{({R^\dag })^2} - 4\mathbf{1}^\top G\mathbf{1}\mathbb{E}[{\lambda^\top  }(B(\omega ) - {R^\dag }r_0(\omega))]\\
	~~~~- {(b({R^\dag }) + 2\mathbf{1}^\top G\mathbf{1}R_\infty ^{ - 1}({R^\dag }))^2}\Bigg],
	\end{array}
	\end{equation}
	where $R^\dagger$ is characterized by
	\begin{equation}	-\frac{dR_\infty^{-1}(R^\dagger)}{dR}=\frac{R_\infty^{-1}(R^\dagger)\mathbf{1}^\top\bar{r}_0-\mathbb{E}[\lambda^\top r_0(\omega)]}{b(R^\dagger)+2\mathbf{1}^\top G\mathbf{1}R_\infty^{-1}(R^\dagger)},
	\end{equation}
	$b(R)=-\bar{\lambda}^\top G\mathbf{1}-\mathbf{1}^\top(\mathbb{E}[B(\omega)]-R\bar{r}_0)$ and $\bar{r_0}$ the expected renewable generation per unit-capacity installed.
\end{corollary}

\begin{proof}[Proof of Corollary \ref{cor:deathflat}]

	We look for the infimum of retailer costs that induce a death spiral. With Proposition \ref{prop:converge} and \ref{prop:piconvex}, such $\theta^\dagger$ is specified when the market potential function $p(R,\theta)$ is tangent to $p=R$, or when  $\pi_{\rm  F}^*(R)$ tangent to $R_\infty^{-1}(R)$. Thus the tangent point can be specified by
	\begin{equation}
	\left\{ \begin{array}{l}
	\pi_{\rm  F}^*(R)^\prime  - {R_\infty^{-1}(R)^{\prime }} = 0\\
	\pi_{\rm  F}^*(R) - {R_\infty^{-1}}(R) = 0
	\end{array} \right.
	\label{eq:tgtpt}
	\end{equation}
	Further deduction of the first equation yields
	\begin{equation}
	\begin{array}{l}
	\pi_{\rm  F} ^*{(R)^\prime } - R_\infty ^{ - 1}{(R)^\prime } = \frac{1}{{\sqrt {{b^2} - 4ac} }}(\pi_{\rm  F} ^*(R){b^\prime } + {c^\prime }) - R_\infty ^{ - 1}{(R)^\prime }\\
	= \frac{1}{{\sqrt {{b^2} - 4ac} }}(R_\infty ^{ - 1}(R){b^\prime } + {c^\prime }) - R_\infty ^{ - 1}{(R)^\prime }
	\end{array}
	\label{eq:fgprime}
	\end{equation}
	Reformulate (\ref{eq:fgprime}) as
	\begin{equation}
	\theta=\frac{1}{4a}[b^2-4ac_0-(\frac{R_\infty ^{ - 1}{(R) }b^\prime+c^\prime}{R_\infty ^{ - 1}{(R)^\prime }})^2]
	\label{eq:f1}
	\end{equation}
	where $c_0=\mathbb{E}[\lambda^\top (B(\omega)-R r_0(\omega))]$. Reformulating the second equation in (\ref{eq:tgtpt}) yields
	\begin{equation}
	\theta=\frac{1}{4a}[b^2-4ac_0-(b+2aR_\infty ^{ - 1}{(R) })^2]
	\label{eq:f2}
	\end{equation}
	With (\ref{eq:f1}) and (\ref{eq:f2}), we can solve $R^\dagger$ which is characterized by
	\begin{equation}
	-\frac{dR_\infty ^{ - 1}(R^\dagger)}{dR}=\frac{R_\infty ^{ - 1}(R^\dagger)\mathbf{1}^\top\bar{r}_0-\mathbb{E}[\lambda^\top r_0(\omega)]}{b+2aR_\infty ^{ - 1}(R^\dagger)}
	\end{equation}
	Substituting $R^\dagger$ into (\ref{eq:f2}) we have
	\begin{equation}
	\begin{array}{l}
	{\theta ^\dag } = \frac{1}{{4a}}[b{({R^\dag })^2} - 4a\mathbb{E}[{\lambda^\top  }(B(\omega ) - {R^\dag }r_0(\omega))]\\
	- {(b({R^\dag }) + 2aR_\infty ^{ - 1}({R^\dag }))^2}],
	\end{array}
	\end{equation}

\end{proof}

\begin{proof}[Proof of Theorem \ref{th:stability}]
	$f(\sigma^*,\chi)=\sigma^*$ (the equilibrium condition) directly holds by computing the dynamics in (\ref{TfromR}) and (\ref{Rdynamics}). 
	
	We prove the stability by constructing the following Lyapunov candidate function. Let $\sigma=(T,R)$,
	\begin{equation}
	V(\sigma)=\max (R^*,R,p(R))-\min (R,R^*).
	\end{equation}
	
	A neighborhood $\mathscr{B}$ of the equilibrium $(T^*,R^*)$ with a range of solar capacity $(v_1,v_2)$ is selected such that (i) $v_2\leq R^*+\epsilon$ and (ii) $p(R)$ is monotonic with $p(R)<v_2$ in $R\in (v_1,R^*)$. We check the left and right regions of $R^*$ respectively:
	
	\textbf{Left region} $R\in (v_1,R^*)$: for a state in the left region, there are four possible situations,
	
	a) if $p(R)\leq R$,
	\begin{equation}
	V(\sigma)=R^*-R>0
	\end{equation}
	\begin{equation}
	V(f(\sigma))-V(\sigma)=R^*-h(R,\mu)-R^*+R= 0
	\end{equation}
	
	b) if $R<p(R)\leq R^*$
	\begin{equation}
	V(\sigma)=R^*-R>0
	\end{equation}
	\begin{equation}
	V(f(\sigma))-V(\sigma)=R^*-h(R,\mu)-R^*+R<0
	\end{equation}
	
	c) if $p(R)> R^*$ and $h(R,\mu)<R^*$
	\begin{equation}
	V(\sigma)=p(R)-R>0
	\end{equation}
	\begin{equation}
	V(f(\sigma))-V(\sigma)=p(h(R,\mu))-h(R,\mu)-p(R)+R<0
	\end{equation}
	
	d) if $p(R)> R^*$ and $h(R,\mu)\geq R^*$
	\begin{equation}
	V(\sigma)=p(R)-R>0
	\end{equation}
	\begin{equation}
	V(f(\sigma))-V(\sigma)=h(R,\mu)-R^*-p(R)+R<0
	\end{equation}

	\textbf{Right region} $R\in (R^*,v_2)$:
	\begin{equation}
	V(\sigma)=R-R^*>0
	\end{equation}
	\begin{equation}
	V(f(\sigma))-V(\sigma)=h(R,\mu)-R^*-R+R^*= 0
	\end{equation}
	
	Thus $(T^*,R^*)$ is a Lyapunov stable equilibrium.
	
	Leveraging Proposition \ref{prop:converge}, it is clear that $\mathop {\lim }\limits_{k \to \infty } {\sigma _k} = {\sigma ^*}$. Since $\sigma^*$ is a stable equilibrium, the adoption is a stable adoption.
		
\end{proof}

\begin{proof}[Proof of Theorem \ref{th:ramseytwopart}]
	
	For Ramsey two-part tariff, the solution of (\ref{eq:ramplan}) has the following expression for volumetric charge (we give the flat and dynamic expressions respectively)
	\begin{equation}
	\label{eq:twopartflat}
	\pi_{\rm F}^\dagger(R,\theta)=\frac{{\mathbb{E}[{\lambda ^ \top }\partial D({\bf{1}}\pi _{\rm{F}}^*,\omega )/\partial \pi _{\rm{F}}^*]}}{\mathbf{1}^\top \mathbb{E} [\partial D(\mathbf{1}\pi_{\rm F}^\dagger,\omega)/\partial \pi_{\rm F}^\dagger]}
	\end{equation}
	\begin{equation}
	\label{eq:twopartdyn}
	\pi^\dagger(R,\theta)=\mathbb{E}{[{\nabla _\pi }D({\pi ^\dag },\omega )]^{ - 1}}\mathbb{E}[{\nabla _\pi }D({\pi ^\dag },\omega )\lambda ]
	\end{equation} 
	Expression (\ref{eq:twopartflat} - \ref{eq:twopartdyn}) reveals that the volumetric rate of Ramsey two-part tariff only depends on the wholesale market prices and the demand function, thus stays unchanged with renewable adoption.
 The market potential function $p(R)$ thus also has the same value for different $R$, that said, is horizontal. The market potential at the equilibrium is then equal to $p(0)$, the market potential at $R=0$. 
 Utilizing Theorem \ref{th:stability}, this equilibrium is stable.  
\end{proof}

\begin{proof}[Proof of Theorem \ref{deflim}]
	We prove the theorem by showing 1) such a $\theta^\circ$ that satisfies $p(R_{\mu_{\rm L}}^\sharp(\theta^\circ))=R_{\mu_{\rm L}}^\sharp(\theta^\circ)$ exists  2) for varying $\theta$,  the highest solution of $p(R)=R$ is $R^\circ$; 3) $R^\circ$ is approachable by a stable adoption process.
	
	1): It is not hard to observe that as $\theta$ increases, $R_{\mu_{\rm L}}^\sharp(\theta)$ approaches to 0. Since $p(R_{\mu_{\rm L}}^\sharp(\theta)$ is bounded by the number of total consumers, leveraging Lemma \ref{lem:criticalconnection}, there exists $\bar{\theta}$ such that $ \forall \theta\geq 0, R_{\mu_{\rm L}}^\sharp(\bar{\theta})\geq p(R_{\mu_{\rm L}}^\sharp(\theta))$. By Brouwer's fixed-point theorem, there exists $\theta^\circ$ that satisfies $p(R_{\mu_{\rm L}}^\sharp(\theta^\circ))=R_{\mu_{\rm L}}^\sharp(\theta^\circ)$.
	
	2): Suppose there exist $R_0>R^\circ$ and $\theta_0$ such that $p(R_0,\theta_0)=R_0$. According to the second assumption of Theorem 4, there must exist $R_{\mu_{\rm L}}^\sharp(\theta_0) \geq R_0$ with $p(R_{\mu_{\rm L}}^\sharp(\theta_0))\geq p(R_0,\theta_0)$. Thus we have $R_{\mu_{\rm L}}^\sharp(\theta_0) \geq R_0>R^\circ$ and $p(R_{\mu_{\rm L}}^\sharp(\theta_0))\geq p(R_0,\theta_0)=R_0>R^\circ=p(R^\circ,\theta^\circ)$, which contradict with the first assumption of Theorem \ref{deflim} and Proposition \ref{prop:rsharpdecrea}.
	
	3): When increasing $\theta$ starting from $\theta^\circ$, leveraging the first assumption in Theorem \ref{deflim} and Proposition \ref{prop:rsharpdecrea}, $p(R_{\mu_{\rm L}}^\sharp(\theta))$ increases while $R_{\mu_{\rm L}}^\sharp(\theta)$ decreases. With the second assumption of Theorem \ref{deflim}, there exists a tariff policy (connection charge $A$ can be varied) such that  $\exists \epsilon>0,R<p(R)\leq R^\circ$ for all $R\in (R^\circ-\epsilon,R^\circ)$. By Theorem \ref{th:stability}, $R^\circ$ is the solar capacity associated with the stable equilibrium of a stable adoption under the tariff policy. 
	
\end{proof}

\begin{proof}[Proof of Lemma \ref{lem:criticalconnection}]	
	Denote 
	\begin{equation}
	\overline {\rm rs}^{\rm M}(A,\theta,R)=\mathop {\max }\limits_\pi \overline {\rm rs} (T,\theta,R).
	\end{equation}
	 $\theta^\sharp$ can be represented by $\theta^\sharp=\mathop {\min }\limits_R  (\overline {\rm rs}^{\rm M}(A,\theta,R)-MA+\theta)$. Let first consider linear tariff, i.e., the connection charge $A=0$. We have $\theta^\sharp=\mathop {\min }\limits_R  (\overline {\rm rs}^{\rm M}(0,\theta,R)+\theta)$. Expanding the $\min$ operation yields
	 \begin{equation}
	  \overline {\rm rs}^{\rm M}(0,\theta,R)+\theta\geq\theta^\sharp, \forall R,
	 \end{equation}
	 which is equivalent to 
	 \begin{equation}
	 \overline {\rm rs}^{\rm M}(0,\theta^\sharp,R)\geq 0, \forall R.
	 \end{equation}
	 With Definition \ref{def:cridif}, the critical adoption level $R^\sharp$ does not exist for linear tariffs under break-even condition. The nonexistence of $R^\sharp$ means there is no death spiral. And since the solar capacity is nondecreasing and the market potential is bounded, there must exist a stable adoption.
	 This result can be easily extended to the two part tariff with fixed connection charge $A^\sharp= (\theta-\theta^\sharp)/M$ by (\ref{eq:rsexpression}).
	
\end{proof}

\end{document}